\documentclass[12pt,a4paper,reqno]{amsart}  
\usepackage{amssymb} \usepackage{amscd} \usepackage{times}
\newcommand{\be} {\begin{equation}}
\newcommand{\ee} {\end{equation}}
\newcommand{\bea} {\begin{eqnarray}}
\newcommand{\eea} {\end{eqnarray}}
\newcommand{\Bea} {\begin{eqnarray*}}
\newcommand{\Eea} {\end{eqnarray*}}

\def\zbb{\mathbb{Z}}  
  
  \def\phi{\varphi}
 \def\p1{{\mathbb{P}^1_\zbb}}

\newtheorem{Theorem}{\quad Theorem}[section]

\newtheorem{Corollary}[Theorem]{\quad Corollary}
\newtheorem{Lemma}[Theorem]{\quad Lemma}
\begin{document}
\title{ Harnack type inequality on Riemannian manifolds of dimension 5.}
\author{Samy Skander Bahoura}
\address{Departement de Mathematiques, Universite Pierre et Marie Curie, 2 place Jussieu, 75005, Paris, France.}
\email{samybahoura@yahoo.fr}
\date{\today}
\maketitle
\begin{abstract}
We give an estimate of type $ \sup \times \inf $ on Riemannian manifold of dimension 5 for the Yamabe equation.
\end{abstract}
%{\bf \small  Keywords: $ \sup \times  \inf $ , riemannian manifold, dimension 3, prescribed curvature, uniqueness.}
\section{Introduction and Main Results}
In this paper, we deal with the following Yamabe equation in dimension $ n= 5 $:
 \be -\Delta_g u+ \dfrac{n-2}{4(n-1)} R_g u=n(n-2) u^{N-1},\,\, u >0, \,\, {\rm and } \,\, N=\dfrac{n+2}{n-2}.  \label{(1)}\ee

Here, $ R_g $ is the scalar curvature.

The equation $ (\ref{(1)}) $ was studied a lot, when $ M =\Omega \subset {\mathbb R}^n $ or $ M={\mathbb S}_n $ see for example, [2-4], [11], [15]. In this case we have a $ \sup \times \inf $ inequality.
The corresponding equation in two dimensions on open set $ \Omega $ of $ {\mathbb R}^2 $, is:
 \be -\Delta u=V(x)e^u, \label{(2)} \ee
The equation $ (\ref{(2)}) $ was studied by many authors and we can find
very important result about a priori estimates in [8], [9], [12],
[16], and [19]. In particular in [9] we have the following interior
estimate:
 $$  \sup_K u  \leq c=c(\inf_{\Omega} V, ||V||_{L^{\infty}(\Omega)}, \inf_{\Omega} u, K, \Omega).  $$
And, precisely, in [8], [12], [16], and [20], we have:
 $$ C \sup_K u + \inf_{\Omega} u \leq c=c(\inf_{\Omega} V, ||V||_{L^{\infty}(\Omega)}, K, \Omega),  $$
and,
 $$ \sup_K u + \inf_{\Omega} u \leq c=c(\inf_{\Omega} V, ||V||_{C^{\alpha}(\Omega)}, K, \Omega).  $$
where $ K $ is a compact subset of $ \Omega $, $ C  $ is a positive constant which depends on $\dfrac{\inf_{\Omega} V}{\sup_{\Omega} V} $, and,  $ \alpha \in (0, 1] $.
When $ \dfrac{4(n-1)h}{n-2}= R_g $ the scalar curvature, and $ M $ compact, the equation $(\ref{(1)}) $ is Yamabe equation. T. Aubin and R. Schoen have proved the existence of solution in this case, see for example [1] and [14] for a complete and detailed summary.
When $ M $ is a compact Riemannian manifold, there exist some compactness result for equation  $ (\ref{(1)}) $ see [18]. Li and Zhu see [18], proved that the energy is bounded and if we suppose $ M $ not diffeormorfic to the three sphere, the solutions are uniformly bounded. To have this result they use the positive mass theorem.
Now, if we suppose $ M $ Riemannian manifold (not necessarily compact) %$ V\equiv 1 $,
 Li and Zhang [17] proved that the product $ \sup \times \inf $ is bounded. Here we extend the result of [5].
 Our proof is an extension Li-Zhang result in dimension 3, see [3] and [17], and,  the moving-plane method is used to have this estimate. We refer to Gidas-Ni-Nirenberg for the  moving-plane method, see  [13]. Also, we can see in [3, 6, 11, 16, 17, 10], some applications of this method, for example an uniqueness result. We refer to [7] for the uniqueness result on the sphere and in dimension 3.
Here, we give an equality of type $ \sup \times \inf $ for the equation $ (\ref{(1)}) $ in dimension 5.
In dimension greater than 3 we have other type of estimates by using moving-plane method, see for example [3, 5].
There are other estimates of type $ \sup + \inf $ on complex Monge-Ampere equation on compact  manifolds, see [20-21] . They consider, on compact Kahler manifold $ (M, g) $, the following equation:
 \be 
 \begin{cases}
(\omega_g+\partial \bar \partial \phi)^n=e^{f-t\phi} \omega_g^n, \\
\omega_g+\partial \bar \partial \phi >0 \,\, {\rm on } \,\,  M  \\
\end{cases}  
  \ee

And, they prove some estimates of type $ \sup_M+m \inf_M \leq C $ or $  \sup_M + m \inf_M \geq C $ under the positivity of the first Chern class of M.
Here, we have,

\begin{Theorem} For all compact set $ K $ of $ M $, there is a positive constant c, which depends only on, $ K, M, g $ such that:
$$ (\sup_K  u)^{1/3}  \times \inf_M u \leq c,  $$
for all $ u $ solution of $ (\ref{(1)}) $.
\end{Theorem}
This theorem extend to the dimension 5 the result of Li and Zhang, see
[17] . Here, we use the method of Li and Zhang in [17] . Also, we extend a result of [5].

\begin{Corollary}  For all compact set $ K $ of $ M $ there is a positive constant c, such that:

 $$ \sup_K  u \leq c=c(g, m, K, M) \,\, {\rm if } \,\, \inf_M u \geq m >0,  $$

for all $ u $ solution of $ (\ref{(1)}) $.
\end{Corollary}

\section{Proof of the theorems}
\underbar {\it Proof of theorem 1.1:}
 We want to prove that
\be \epsilon^{3} (\max_{B(0, \epsilon)} u)^{1/3} \times \min_{B(0, 4
\epsilon) } u \leq c=c(M, g) \label{(3)}.\ee We argue by
contradiction and we assume that \be\label{4} ( \max_{B(0,
\epsilon_k)} u_k)^{1/3} \times \min_{B(0, 4 \epsilon_k) } u_k \geq k
{\epsilon_k }^{-3}.\ee \underbar {Step 1: The blow-up analysis }
The blow-up analysis gives us : For some $  \bar x_k \in B(0,
\epsilon_k) $, $ u_k( \bar x_k)=\max_{B(0, \epsilon_k)} u_k $, and,
from the hypothesis,
$$ u_k( \bar x_k)^{4/9} \epsilon_k \to + \infty. $$
By a standard selection process, we can find $ x_k  \in B(\bar x_k ,
\epsilon_k/2) $ and $ \sigma_k \in (0, \epsilon_k/4) $ satisfying,
\bea
 u_k(x_k)^{4/9}\sigma_k &\to &+ \infty, \label{(5)}\\
  u_k(x_k) &\geq &  u_k(\bar x_k), \label{(6)}\\
 {\rm and},~~ u_k(x) &\leq &  C u_k (x_k),  \,\, {\rm in} \,\, B(x_k, \sigma_k),
\label{(7)}\eea
 where $ C $ is some universal constant. It follows
from above (\ref{4}), (\ref{(6)}) that \be
( u_k(x_k))^{1/3} \times
(\min_{\partial B(x_k, 2 \epsilon_k) } u_k) \sigma_k^3  \geq  (u_k(\bar
x_k))^{1/3} \times (\min_{B(0, 4 \epsilon_k) } u_k)\epsilon_k^3 \geq k   \to +
\infty .\label{(8)}\ee 

We use $ \{ z^1, \ldots, z^n \}  $ to denote
some geodesic normal coordinates centered at $ x_k $ (we use the
exponential map). In the geodesic normal coordinates, $
g=g_{ij}(z)dz^idz^j $, \be g_{ij}(z)-{\delta}_{ij} = O(r^2), \,\,
g:=det(g_{ij}(z))=1+O(r^2), \,\, h(z)=O(1), \label{(9)}\ee

where $
r= |z| $. Thus,
$$ \Delta_g u=\dfrac{1}{\sqrt g} \partial_i ( \sqrt g g^{ij} \partial_j u)=\Delta u+ b_i\partial_i u +  d_{ij} \partial_{ij} u, $$
where \be b_j=O(r),  \,\, d_{ij}= O(r^2) \label{(10)}\ee We have a
new function
$$ v_k(y)=M_k^{-1} u_k(M_k^{-2/(n-2)}y) \,\, {\rm for} \,\, |y| \leq 3 \epsilon_k  M_k^{2/(n-2)} $$
where $ M_k=u_k(0) $. From (\ref{(5)}) and (\ref{(8)}) we have \be
\left.
\begin{array}{lll}
 \Delta v_k + \bar b_i\partial_i v_k + \bar d_{ij} \partial_{ij} v_k - \bar c v_k +{v_k}^{N-1} &= &0  \,\, {\rm for} \,\, |y| \leq 3 \epsilon_k M_k^{2/(n-2)}  \\
 v_k(0)&=& 1  \\
  v_k(y)& \leq & C_1 \, {\rm for} \,\, |y| \leq  \sigma_k  M_k^{2/(n-2)}\\
   \lim_{k \to + \infty } \min_{|y|=2\epsilon_k M_k^{4/9}} (v_k(y)|y|^3 )= + \infty. \end{array} \right\}  \label{(11)}\\ \ee
where $ C_1 $ is a universal constant and  \be
 \bar b_i(y)=
M_k^{-2/(n-2)} b_i(M_k^{-2/(n-2)}y), \,\, \bar
d_{ij}(y)=d_{ij}(M_k^{-2/(n-2)}y) \label{(12)}  \ee and,  \be
 \bar
c(y)=M_k^{-4/(n-2)} h(M_k^{-2/(n-2)}y) \label{(13)}.\ee
 We can see
that for $ |y| \leq 3 \epsilon_k M_k^{2/(n-2)}  $, \be |\bar b_i(y)|
\leq CM_k^{-4/(n-2)}|y|, \,\, |\bar  d_{ij}(y)| \leq CM_k^{-4/(n-2)}
|y|^2, \,\, |\bar c(y)| \leq CM_k^{-4/(n-2)} \label{(14)}\ee where $
C $ depends on $ n, M, g $.

It follows from (\ref{(11)}), (\ref{(12)}), (\ref{(13)}), (\ref{(14)})
and the elliptic estimates, that, along a subsequence, $ v_k $
converges in $ C^2 $ norm on any compact subset of $ {\mathbb R}^2 $
to a positive function $ U $ satisfying
 \be
\left. \begin{array}{lll}
\Delta U+ U^{N-1} & = & 0,\,\,  {\rm in } \,\,  {\mathbb R}^n, \,\,{\rm with } \,\, N=\dfrac{n+2}{n-2} \\
U(0)=1, &&  0 < U \leq  C.   \end{array}\right\} \label{(15)} 
 \ee
 
 In the case where $ C=1 $, by a result of Caffarelli-Gidas-Spruck, see [10], we have:
 
 \be U(y)=(1+|y|^2)^{-(n-2)/2}, \ee

But, here we do not need this result.

Now, we need a precision in the previous estimates, we take a conformal change of metric such that, the Ricci tensor vanish,

 \be R_{jp} = 0. \ee

We have by the expressions for $ g $ and $  g_{ij} $, as in the paper of Li-Zhang,

 \be  b_j=O(r^2),   \,\, R= O(r),  \,\, d_{ij}=-\dfrac{1}{3}R_{ipqj}z^p z^q + O(r^3). \ee 
 
 Thus,
 
 \be |\bar c| \leq C |y|M_k^{-2}, \,\, |\bar b_i| \leq C|y|^2M_k^{-2}, \ee 
 
and,

\be \bar d_{ij}=-\dfrac{1}{3}M_k^{-4/3}R_{ipqj}y^p y^q + O(1)M_k^{-2} |y|^{3}. \ee

As, in the paper of Li-Zhang, we have:

\be v_k(y)\geq C|y|^{-3},  \,\, 1 \leq |y| \leq 2 \epsilon_k M_k^{2/3}. \label{(16)}  \ee

with $ C>0 $.

For $ x\in {\mathbb R}^2 $ and $ \lambda >0 $, let,

\be v_k^{\lambda, x}(y):= \dfrac {\lambda}{|y-x|}v_k \left (x+\dfrac {\lambda^2(y-x)}{|y-x|^2}\right  ), \ee

denote the Kelvin transformation of $ v_k $ with respect to the ball centered at x and of radius $  \lambda  $.

\bigskip

We want to compare for fixed $ x $, $ v_k $ and $ v_k^{\lambda, x} $. For simplicity we assume $ x=0 $. We have:

$$ v_k^{\lambda}(y):= \dfrac {\lambda}{|y|}v_k(y^{\lambda}), \,\, {\rm with } \,\,  y^{\lambda}= \dfrac {\lambda^2y}{|y|^2}. $$

For $ \lambda >0 $, we set,

$$ \Sigma_{\lambda}=B\left (0, \epsilon_k{M_k}^{2}\right )-{\bar B(0,\lambda)}. $$

The boundary condition, $  (\ref{(11)}) $, become:

\be \lim_{k \to + \infty } \min_{|y|=\epsilon_k M_k^{4/9}} \left( v_k(y)|y|^3 \right )= \lim_{k \to + \infty } \min_{|y|=2\epsilon_k M_k^{4/9}} \left ( v_k(y)|y|^3 \right )= + \infty.\ee

As in the paper of Li-Zhang, we have:

\be \Delta w_{\lambda} + \bar b_i\partial_i w_{\lambda} + \bar d_{ij} \partial_{ij} w_{\lambda}- \bar c w_{\lambda} + \dfrac {(n+2)}{(n-2)}\xi^{4/(n-2)} w_{\lambda} = E_{\lambda} \,\, {\rm in  }  \,\, \,\,    \Sigma_{\lambda }.\ee

where $ \xi $ stay between $ v_k $ and $ v_k^{\lambda} $. Here,

$$ E_{\lambda} = -  \bar b_i \partial_i v_k^{\lambda} - \bar d_{ij} \partial_{ij} v_k^{\lambda} + \bar c v_k^{\lambda} -E_1, $$

with,

\be E_1(y)=- \left (\dfrac{\lambda}{|y|} \right )^{n+2} \left (  \bar b_i(y^{\lambda}) \partial_i v_k(y^{\lambda}) + \bar d_{ij}(y^{\lambda}) \partial_{ij} v_k(y^{\lambda}) - \bar c(y^{\lambda}) v_k(y^{\lambda}) \right ). \ee

\begin{Lemma} We have,

\be |E_{\lambda}|\leq C_1 \lambda|y|^{-1}M_k^{-2}+C_2 \lambda^3|y|^{-3}M_k^{-4/3}.\ee
 
\end{Lemma} 

{\bf Proof:} as in the paper of Li-Zhang, we have a nonlinear term $ E_{\lambda} $ with the following property,

$$ |E_{\lambda}|\leq C_1 \lambda^3 M_k^{-2}|y|^{-2}+C_2 \lambda^4 M_k^{-4/3}|y|^{-4}\leq C_1 \lambda|y|^{-1}M_k^{-2}+C_2 \lambda^3|y|^{-3}M_k^{-4/3}. $$

 Next, we need an auxiliary function which correct the nonlinear term. Here we take the following auxiliary function:

\be h_{\lambda}= -C_1\lambda M_k^{-2}\left (|y|-\lambda \right ) - C_2 \lambda^2 M_k^{-4/3}\left (\left (1- \left (\dfrac{\lambda}{|y|}\right )^3\right )-\left (1-\left (\dfrac{\lambda}{|y|}\right )\right )\right ), \ee

we have,

\be h_{\lambda}\leq 0,  \ee

\be \Delta h_{\lambda} = -C_1\lambda|y|^{-1}M_k^{-2}-C_2 \lambda^3|y|^{-3}M_k^{-4/3}, \ee

and, thus,

$$ \Delta h_{\lambda} + |E_{\lambda}| \leq 0. $$

As in the paper of Li-Zhang, we can prove the following lemma:

\begin{Lemma} We have,

\be w_{\lambda } + h_{\lambda } > 0 , \,\, {\rm in  }  \,\, \Sigma_{\lambda }  \,\, \forall 0 < \lambda \leq \lambda_1. \ee
 
\end{Lemma} 

Before to prove the lemma, note that, here, we consider the fact that,

\be \lambda \leq |y|\leq \epsilon_k M_k^{4/9}\leq \epsilon_k M_k^{2/3}. \ee

And, as in the paper of Li-Zhang, we need the estimate $ (\ref{(16)}) $:

$$ v_k(y)\geq C |y|^{-3},  \,\, 1 \leq |y| \leq 2 \epsilon_k M_k^{2/3}. $$

with $ C >0 $.

{\bf Proof : }

{\bf Step 1: } There exists $ \lambda_0 >0 $ independent of $ k $ such the assertion of the lemma holds for all $ 0 <  \lambda < \lambda_0 $.

To see this, we write:

$$ w_{\lambda } = v_k(y)-v_k^{\lambda }(y)=|y|^{-3/2}(|y|^{3/2} v_k(y)-|y^{\lambda }|^{3/2} v_k(y^{\lambda })). $$ 

Let, in polar coordinates,

$$ f (r, \theta)=r^{3/2}v_k(r, \theta). $$

By the properties of $ v_k $, there exist $ r_0 >0 $ and $ C>0 $ independant of $ k $ such that:

$$ \partial_r f (r, \theta) > Cr^{1/2},  \,\, {\rm for  }  \,\, 0 <r < r_0. $$

Thus,

$$ w_{\lambda }(y)  =|y|^{-3/2} (f(|y|, y/|y|)-f(|y^{\lambda }|,y/|y|))=|y|^{-3/2} \int_{|y^{\lambda }|}^{|y|} \partial_r f (r, y/|y|) dr >  $$

$$ > C' |y|^{-3/2}(|y|^{3/2} - |y^{\lambda }|^{3/2}) > C'' (|y| - \lambda) \,\, {\rm for  }  \,\, 0< \lambda < |y|< r_0 ,$$ 
    
with, $ C', C'' >0 $.

It follows that,

\be w_{\lambda }+h_{\lambda } \geq  (C''-o(1))(|y| - \lambda), \,\, {\rm for  }  \,\, 0< \lambda < |y|< r_0, \ee

Now, for 

$$ r_0 \leq |y|\leq \epsilon_k M_k^{4/9}\leq \epsilon_k M_k^{2/3}. $$

we have by the definition of $ h_{\lambda} $, and, as in the paper of Li-Zhang, we need the estimate $ (\ref{(16)}) $:

$$ v_k(y)\geq C|y|^{-3},  \,\, 1 \leq |y| \leq 2 \epsilon_k M_k^{2/3}. $$

to have,

$$  |h_{\lambda} | < \dfrac{1}{2} v_k (y). $$

Thus, as in the paper of Li-Zhang,

\be w_{\lambda }+h_{\lambda } >0,  \,\, {\rm for  }  \,\, 0< r_0 < |y|<2 \epsilon_k M_k^{4/9}. \ee

{\bf  Step 2:} Set, 
 
 \be  \bar {\lambda }^k = \sup \{ 0  < \lambda \leq \lambda_1, w_{\mu} + h_{\mu} \geq 0 , \,\, {\rm in  }  \,\, \Sigma_{\mu }  \,\, \forall 0 < \mu \leq \lambda \}, \ee
 
We claim that, $ \bar {\lambda }^k = \lambda_1 $.

In order to apply the maximum principle and the Hopf lemma, we need to prove that:

\be  (\Delta + \bar b_i\partial_i  + \bar d_{ij} \partial_{ij} - \bar c) ( w_{\lambda} + h_{\lambda })\leq 0 \,\, {\rm in  }  \,\, \,\,    \Sigma_{\lambda } \label{(17)}\ee

In other words, we need to prove that:

\be \Delta h_{\lambda} + \bar b_i\partial_i h_{\lambda} + \bar d_{ij} \partial_{ij} h_{\lambda}- \bar c h_{\lambda}+ E_{\lambda }\leq 0 \,\, {\rm in  }  \,\, \,\,    \Sigma_{\lambda } .\ee

First note that, $ h_{\lambda} < 0 $. Here, we consider the fact that,

$$ \lambda \leq |y|\leq \epsilon_k M_k^{4/9}\leq \epsilon_k M_k^{2/3}. $$

We have,

$$ |\bar c| \leq C |y|M_k^{-2}, $$

Thus,

$$ |y||\bar c h_{\lambda}| \leq C_1M_k^{-4}\lambda |y|^2(|y|-\lambda)+C_2M_k^{-10/3}\lambda |y|^2 \leq o(1)M_k^{-2}\lambda, $$

which we can write as,

\be |\bar c h_{\lambda}| \leq C_1M_k^{-2}\lambda |y|^{-1}. \ee 

We have,

$$  |\bar b_i| \leq C|y|^2M_k^{-2}, $$

Thus,

$$ |\bar b_i  \partial_i h_{\lambda}| \leq C_1M_k^{-4}\lambda |y|^2+C_2M_k^{-10/3}(\lambda^5 |y|^{-2}+ \lambda^3),   $$

$$ |y| C_1M_k^{-4}\lambda |y|^2=o(1)M_k^{-2}\lambda, $$

which we can write as,

\be C_1M_k^{-4}\lambda |y|^2=o(1)M_k^{-2}\lambda |y|^{-1}. \ee 

and,

$$ |y|^3 C_2M_k^{-10/3}\lambda^5|y|^{-2}=C_2M_k^{-10/3}\lambda^3 |y|= o(1)M_k^{-4/3}\lambda^3, $$

which we can write as,

 \be C_2M_k^{-10/3}\lambda^5|y|^{-2} = o(1)M_k^{-4/3}\lambda^3 |y|^{-3}. \ee

and,

$$ |y|^3 C_2M_k^{-10/3}\lambda^3=o(1)M_k^{-2}\lambda^3, $$

which we can write as,

\be C_2M_k^{-10/3}\lambda^3 = o(1)M_k^{-4/3}\lambda^3 |y|^{-3}. \ee

Thus,

\be |\bar b_i  \partial_i h_{\lambda}| \leq  o(1)M_k^{-2}\lambda |y|^{-1}+o(1)M_k^{-4/3}\lambda^3 |y|^{-3}.  \ee

We have,

$$  |\bar d_{ij}| \leq |y|^2M_k^{-4/3}, $$

Thus,

$$  |\bar d_{ij} \partial_{ij} h_{\lambda}| \leq \lambda |y| M_k^{-10/3} +C_2M_k^{-8/3}(\lambda^5 |y|^{-3}+ \lambda^3|y|^{-1}) ,  $$

Thus,

$$  |\bar d_{ij} \partial_{ij} h_{\lambda}| \leq \lambda |y|^{-1} M_k^{-10/3} +o(1)M_k^{-8/3}\lambda^5 |y|^{-3}+o(1)M_k^{-8/3} \lambda^3 |y|^{-1} ,  $$

Finaly,

\be  |\bar d_{ij} \partial_{ij} h_{\lambda}| \leq o(1)\lambda |y|^{-1} M_k^{-2} +o(1)M_k^{-4/3}\lambda^3 |y|^{-3}+o(1)M_k^{-2} \lambda |y|^{-1}. \ee

Finaly,

\be  |\bar d_{ij} \partial_{ij} h_{\lambda}+\bar b_i  \partial_i h_{\lambda}+\bar c h_{\lambda}|\leq o(1)\lambda |y|^{-1} M_k^{-2} +o(1)\lambda^3 |y|^{-3}M_k^{-4/3}. \ee

Finaly, we have,

$$ \Delta h_{\lambda} + \bar b_i\partial_i h_{\lambda} + \bar d_{ij} \partial_{ij} h_{\lambda}- \bar c h_{\lambda}+ E_{\lambda }\leq 0 \,\, {\rm in  }  \,\, \,\,    \Sigma_{\lambda }, $$

And, thus $ (\ref{(17)}) $,

$$  (\Delta + \bar b_i\partial_i  + \bar d_{ij} \partial_{ij} - \bar c) ( w_{\lambda} + h_{\lambda })\leq 0 \,\, {\rm in  }  \,\, \,\,    \Sigma_{\lambda }. $$

Also, we have from the boundary condition and the definition of $ v_k^{\lambda } $ and $ h_{\lambda} $, we have:

\be  | h_{\lambda } (y)| +v_k^{\lambda }(y) \leq  \dfrac{ C(\lambda_1) }{|y|^3 } ,   \,\,\,\, \forall   \,\,  |y|= \epsilon_k M_k^{4/9},    \ee

and, thus,

\be  w_{\bar \lambda^k }(y) + h_{ \bar \lambda^k } (y) >0 \,\,\,\, \forall   \,\,  |y|= \epsilon_k M_k^{4/9},    \ee

We can use the maximum principle and the Hopf lemma to have:

\be w_{\bar \lambda^k } + h_{ \bar \lambda^k } >0, \,\, {\rm in  }  \,\, \,\, \Sigma_{\lambda } , \ee
and,

\be  \dfrac{\partial }{ \partial \nu} (w_{\bar \lambda^k } + h_{ \bar \lambda^k }) >0, \,\, {\rm in  }  \,\, \,\, \Sigma_{\lambda }. \ee

From the previous estimates we conclude that $ \bar \lambda^k = \lambda_1 $ and the lemma is proved.

\bigskip

Given any  $ \lambda >0 $, since  the sequence $ v_k $ converges to $ U $ and $ h_{ \bar \lambda^k } $ converges to 0 on any compact subset of $ {\mathbb R}^2 $, we have:

\be U(y) \geq U^{\lambda}(y), \,\,\, \forall   \,\,  |y| \geq \lambda, \,\, \forall   \,\,  0 < \lambda  < \lambda_1. \ee

Since $ \lambda_1>0 $ is arbitrary, and since we can apply the same argument to compare $ v_k $ and $ v_k^{\lambda, x} $,  we have:

\be U(y) \geq U^{\lambda, x}(y), \,\,\, \forall   \,\,  |y-x| \geq \lambda >0. \ee

Thus implies that $ U $ is a constant which is a contradiction.


\bigskip

\begin{thebibliography}{99}
\bibitem{1}{T. Aubin. Some Nonlinear Problems in Riemannian Geometry. Springer-Verlag 1998 }
\bibitem{2}{ S.S Bahoura. Majorations du type $ \sup u \times \inf u \leq c $ pour l'\'equation de la courbure scalaire sur un ouvert de $ {\mathbb R}^n, n\geq 3 $. J. Math. Pures. Appl.(9) 83 2004 no, 9, 1109-1150.}
\bibitem{3}{S.S. Bahoura. Harnack inequalities for Yamabe type equations.  Bull. Sci. Math.  133  (2009),  no. 8, 875-892}
\bibitem{4}{S.S. Bahoura. Lower bounds for sup+inf and sup $ \times $ inf and an extension of Chen-Lin result in dimension 3.  Acta Math. Sci. Ser. B Engl. Ed.  28  (2008),  no. 4, 749-758}
\bibitem{5}{S.S. Bahoura. Estimations uniformes pour l'equation de Yamabe en dimensions 5 et 6. J. Funct. Anal.  242  (2007),  no. 2, 550-562.}
\bibitem{6}{S.S. Bahoura. sup $ \times $ inf  inequality on manifold of dimension 3, to appear in MATHEMATICA AETERNA}
\bibitem{7}{H. Brezis, YY. Li. Some nonlinear elliptic equations have only constant solutions. J. Partial Differential Equations 19 (2006), no. 3, 208-217.
}
\bibitem{8}{H. Brezis, YY. Li , I. Shafrir. A sup+inf inequality for some
nonlinear elliptic equations involving exponential
nonlinearities. J.Funct.Anal.115 (1993) 344-358.
}
\bibitem{9}{H.Brezis and F.Merle, Uniform estimates and blow-up behavior for solutions of $ -\Delta u=Ve^u $ in two dimensions, Commun Partial Differential Equations 16 (1991), 1223-1253.
}
\bibitem{10}{L. Caffarelli, B. Gidas, J. Spruck. Asymptotic symmetry and local
behavior of semilinear elliptic equations with critical Sobolev
growth. Comm. Pure Appl. Math. 37 (1984) 369-402.
}
\bibitem{11}{C-C.Chen, C-S. Lin. Estimates of the conformal scalar curvature
equation via the method of moving planes. Comm. Pure
Appl. Math. L(1997) 0971-1017.}
\bibitem{12}{C-C.Chen, C-S. Lin. A sharp sup+inf inequality for a nonlinear elliptic equation in ${\mathbb R}^2$.
Commun. Anal. Geom. 6, No.1, 1-19 (1998).}
\bibitem{13}{B. Gidas, W-Y. Ni, L. Nirenberg. Symmetry and related properties via the maximum principle.  Comm. Math. Phys.  68  (1979), no. 3, 209-243.}
\bibitem{14}{J.M. Lee, T.H. Parker. The Yamabe problem. Bull.Amer.Math.Soc (N.S) 17 (1987), no.1, 37 -91.}
\bibitem{15}{YY. Li. Prescribing scalar curvature on $ {\mathbb S}_n $ and related
Problems. C.R. Acad. Sci. Paris 317 (1993) 159-164. Part
I: J. Differ. Equations 120 (1995) 319-410. Part II: Existence and
compactness. Comm. Pure Appl.Math.49 (1996) 541-597.
}
\bibitem{16}{YY. Li. Harnack Type Inequality: the Method of Moving Planes. Commun. Math. Phys. 200,421-444 (1999).}
\bibitem{17}{YY. Li, L. Zhang. A Harnack type inequality for the Yamabe equation in low dimensions.  Calc. Var. Partial Differential Equations  20  (2004),  no. 2, 133--151.}
\bibitem{18}{YY.Li, M. Zhu. Yamabe Type Equations On Three Dimensional Riemannian Manifolds. Commun.Contem.Mathematics, vol 1. No.1 (1999) 1-50.}
\bibitem{19}{I. Shafrir. A sup+inf inequality for the equation $ -\Delta u=Ve^u $. C. R. Acad.Sci. Paris S\'er. I Math. 315 (1992), no. 2, 159-164.}
\bibitem{20}{Y-T. Siu. The existence of Kahler-Einstein metrics on manifolds with positive anticanonical line bundle and a suitable finite symmetry group.  Ann. of Math. (2)  127  (1988),  no. 3, 585-627}
\bibitem{21}{G. Tian. A Harnack type inequality for certain complex Monge-Ampere equations.  J. Differential Geom.  29  (1989),  no. 3, 481-488.}

\end{thebibliography}
\end{document}